\documentclass[a4paper,11pt]{amsart}
\usepackage[a4paper,twoside,centering]{geometry}
\usepackage{amsfonts}

\title{Hochschild cohomology of gentle algebras}
\author{Sefi Ladkani}
\address{%
Mathematical Institute of the University of Bonn \\
Endenicher Allee 60 \\
53115 Bonn, Germany}
\urladdr{http://www.math.uni-bonn.de/people/sefil}
\email{sefil@math.uni-bonn.de}

\date{August 31, 2011}

\DeclareMathOperator{\ch}{char}
\DeclareMathOperator{\Ext}{Ext}
\DeclareMathOperator{\HH}{HH}

\newcommand{\eps}{\varepsilon}
\newcommand{\gL}{\Lambda}
\newcommand{\bN}{\mathbb{N}}

\theoremstyle{plain}
\newtheorem{theorem}{Theorem}
\newtheorem{cor}{Corollary}

\theoremstyle{definition}
\newtheorem*{remark*}{Remark}
\newtheorem*{notat*}{Notations}

\begin{document}

\begin{abstract}
We compute the Hochschild cohomology groups of gentle algebras and show
that they are determined by the derived invariant introduced by
Avella-Alaminos and Geiss.
\end{abstract}

\maketitle

Gentle algebras are certain finite-dimensional algebras defined
combinatorially in terms of quivers with relations. They have
remarkable homological properties; they are
Gorenstein~\cite{GeissReiten05} and the class of gentle algebras is
closed under derived equivalence~\cite{SchroerZimmermann03}. It is a
long standing problem to classify gentle algebras up to derived
equivalence. To this end, one searches for derived invariants that will
allow to distinguish the different derived equivalence classes.

One such invariant was developed by Avella-Alaminos and
Geiss~\cite{AvellaAlaminosGeiss08}. It takes the form of a function
$\phi_{\gL} : \bN^2 \to \bN$ that can be effectively computed from the
quiver with relations of a gentle algebra $\gL$. In some cases this
invariant is fully capable of distinguishing the derived equivalence
classes~\cite{AvellaAlaminos08,AvellaAlaminosGeiss08,BobinskiMalicki08},
but in general it is not complete and thus further derived invariants
are needed.

A possible candidate is the Hochschild cohomology, which is well-known
to be derived invariant~\cite{Happel89,Rickard91}, defined as
$\HH^{*}(\gL) = \Ext^{*}_{\gL^{op} \otimes_K \gL}(\gL, \gL)$ for an
algebra $\gL$ over a field $K$. Indeed, in~\cite{BessenrodtHolm07} it
was used, together with other derived invariants, to classify up to
derived equivalence the gentle algebras with two vertices and those
with three vertices and zero Cartan determinant.

However, in this note we show that at least as a graded vector space,
the Hochschild cohomology of a gentle algebra is completely determined
by its Avella-Alaminos-Geiss derived invariant and thus cannot
distinguish derived classes not already distinguished by the latter. In
fact, we give an explicit formula for the dimensions of the Hochschild
cohomology groups in terms of the Avella-Alaminos-Geiss derived
invariant, based on the projective resolution of a monomial algebra as
a bimodule over itself given by Bardzell~\cite{Bardzell97} (see
also~\cite{Skoldberg99}), observing that gentle algebras are quadratic
monomial algebras.

In order to formulate our results more precisely, we encode the
dimensions of the Hochschild cohomology groups of a finite-dimensional
algebra $\gL$ over a field $K$ in a formal power series
\[
h_{\gL}(z) = \sum_{i=0}^{\infty} \dim_K \HH^i(\gL) \cdot z^i - 1
\]
and we define, for $n \geq 1$, the formal power series
\[
g_n(z) =
\frac{z^n(1+z)(\eps_n+z^n)}{1-z^{2n}} ,
\quad
\text{where}
\quad
\eps_n = \begin{cases}
1 & \text{if $n$ is even or $\ch K = 2$,} \\
0 & \text{otherwise.}
\end{cases}
\]

Recall that for a finite quiver $Q$, the \emph{Euler characteristic}
$\chi(Q)$ is defined as the number of its vertices minus the number of
its arrows.

\begin{notat*}
For a gentle algebra $\gL$, denote by $Q_{\gL}$ its quiver and by
$\phi_{\gL}$ its Avella-Alaminos-Geiss derived invariant. We assume
that all algebras are \emph{connected}.
\end{notat*}

\begin{theorem}
Let $\gL$ be a gentle algebra. Then
\[
h_{\gL}(z) = \bigl(1 - \chi(Q_{\gL})\bigr) z  + \sum_{m \geq 0}
\phi_{\gL}(1,m) z^m + \sum_{m > 0} \phi_{\gL}(0,m) g_m(z)
\]
\end{theorem}

\begin{remark*}
The Euler characteristic $\chi(Q_{\gL})$ can be expressed in terms of
$\phi_{\gL}$ by the formula
\[
2 \cdot \chi(Q_{\gL}) = \sum_{(n,m) \in \bN^2} \phi_{\gL}(n,m)(n-m) ,
\]
hence $h_{\gL}(z)$ depends only on $\phi_{\gL}$ and whether or not the
field $K$ has characteristic $2$.
\end{remark*}

Let us write explicitly the dimensions of the Hochschild cohomology
groups. We define $\psi_{\gL}(n) = \sum_{d \mid n} \phi_{\gL}(0,d)$ for
$n \geq 1$.

\begin{cor}
Let $\gL$ be a gentle algebra. Then
\begin{enumerate}
\renewcommand{\theenumi}{\alph{enumi}}
\item
$\dim \HH^0(\gL) = 1 + \phi_{\gL}(1,0)$.

\item
$\dim \HH^1(\gL) = 1 - \chi(Q_{\gL}) + \phi_{\gL}(1,1) +
\begin{cases}
\phi_{\gL}(0,1) & \text{if $\ch K = 2$,} \\
0 & \text{otherwise.}
\end{cases}
$

\item
$\dim \HH^n(\gL) = \phi_{\gL}(1,n) +
a_n \psi_{\gL}(n) + b_n \psi_{\gL}(n-1)$ for $n \geq 2$, where
\[
(a_n, b_n) = \begin{cases}
(1,0) & \text{if $\ch K \neq 2$ and $n$ is even,} \\
(0,1) & \text{if $\ch K \neq 2$ and $n$ is odd,} \\
(1,1) & \text{if $\ch K = 2$.}
\end{cases}
\]
\end{enumerate}
\end{cor}

\begin{remark*}
The numbers $\dim \HH^n(\gL)$ can be effectively computed from the
quiver with relations of $\gL$, since the same is true for
$\phi_{\gL}$.
\end{remark*}

\begin{cor}
The following conditions are equivalent for a gentle algebra $\gL$.
\begin{enumerate}
\renewcommand{\theenumi}{\roman{enumi}}
\item
$\HH^1(\gL)=0$;

\item
$\HH^i(\gL)=0$ for all $i \geq 1$;

\item
$\gL$ is piecewise hereditary of Dynkin type $A$.
\end{enumerate}
\end{cor}

\begin{cor}
Let $\gL$ be a gentle algebra of finite global dimension. Then
\[
h_{\gL}(z) = \bigl(1-\chi(Q_{\gL})\bigr)z + \sum_{n \geq 0}
\phi_{\gL}(1,n) z^n .
\]
In particular, $\dim \HH^n(\gL) = \phi_{\gL}(1,n)$ for $n>1$.
\end{cor}

As an application of our results, we compute the Hochschild cohomology
groups of the gentle algebras arising from surface triangulations
introduced by Assem, Br\"{u}stle, Charbonneau and
Plamondon~\cite{ABCP10}.

A marked bordered oriented surface without punctures is a pair $(S,M)$
where $S$ is a compact, connected Riemann surface with non-empty
boundary $\partial S$ and $M \subset \partial S$ is a finite set of
\emph{marked points} containing at least one point from each connected
component of $\partial S$. In~\cite{ABCP10} the authors associate to
each ideal triangulation $T$ of $(S,M)$ a gentle algebra $\gL_T$. The
invariant $\phi_{\gL_T}$ was computed
in~\cite{DavidRoeslerSchiffler11}.

\begin{notat*}
Denote by $g$ the genus of $S$, by $b \geq 1$ the number of connected
components of its boundary $\partial S$ and by $c_1$ the number of
those components with exactly one marked point. For an ideal
triangulation $T$ of $(S,M)$, we set the following quantities:
\begin{enumerate}
\item[$c_0$ --]
the number of boundary components with exactly two marked points such
that their two boundary segments are sides of the same triangle in $T$;

\item[$d$ --]
the number of triangles of $T$ such that two of their sides are
boundary segments.
\end{enumerate}
In addition we set $f_3(z) = z + g_3(z)$ (this agrees with our notation
in~\cite{Ladkani11}).
\end{notat*}

\begin{theorem}
Let $T$ be an ideal triangulation of a marked bordered oriented surface
without punctures. Then
\[
h_{\gL_T}(z) = c_0 + (2g+b-1+c_1)z + \bigl( 4(g-1)+2b+d \bigr) f_3(z) .
\]
\end{theorem}

\begin{cor}
$\HH^2(\gL) = 0$ for any gentle algebra $\gL$ arising from surface
triangulation.
\end{cor}

The above theorem applies in particular to cluster-tilted algebras of
affine type $\widetilde{A}$. In~\cite{Bastian09}, the quivers of these
algebras have been explicitly described and their derived equivalence
classes were characterized in terms of four parameters.

\begin{cor}
Let $\gL$ be a cluster-tilted algebra of type $\widetilde{A}$ with
parameters $(s_1,t_1,s_2,t_2)$. Then
\[
h_{\gL}(z) =  c_0 + (1+c_1)z + (t_1+t_2) f_3(z)
\]
where $c_0 = \left|\{i \,:\, (s_i,t_i)=(0,1)\}\right|$ and $c_1 =
\left|\{i \,:\, (s_i,t_i)=(1,0)\}\right|$.
\end{cor}

\bibliographystyle{amsplain}
\bibliography{hhgentle}

\end{document}